\documentclass[11pt]{article}%
\usepackage{amsmath}
\usepackage{amsfonts}
\usepackage{amssymb}
\usepackage{graphicx}
\usepackage{a4wide}
\providecommand{\U}[1]{\protect\rule{.1in}{.1in}}
\newtheorem{theorem}{Theorem}

\newtheorem{example}[theorem]{Example}

\newtheorem{proposition}[theorem]{Proposition}

\begin{document}

\title{On three open problems related to quasi relative interior}
\author{C. Z\u{a}linescu\thanks{Faculty of Mathematics, University Alexandru Ioan
Cuza, Ia\c{s}i, Romania, e-mail: \texttt{zalinesc@uaic.ro}. } }
\date{}
\maketitle

\begin{abstract}
We give answers to two questions formulated by Borwein and Goebel in
2003 and to a conjecture formulated by Grad and Pop in 2014 related
to calculus rules for quasi (relative) interior.


%\subclass{MSC code1 \and MSC code2 \and more}

\end{abstract}
The notion of quasi relative interior, introduced by Borwein and
Lewis in 1992, became very familiar in the last ten years, being
used mostly for getting necessary optimality conditions in (scalar
or vector) convex programming. Unfortunately the calculus rules for
the quasi relative interior is much poorer than that for other types
of interiority notions. In this short note we give answers to two
questions formulated by Goebel and Borwein in \cite{BorGoe:03} and
to a conjecture of Grad and Pop from \cite{GraPop:14} related to
calculus for quasi (relative) interior.

Throughout the paper, if not specified otherwise, $X$ is a separated
locally convex space and $X^{\ast}$ is its topological dual. For
$x\in X$ and $x^{\ast}\in X^{\ast}$ we set $\left\langle
x,x^{\ast}\right\rangle :=x^{\ast }(x)$. Having $A\subset X$ we use
the notations $\operatorname*{cl}A$, $\operatorname*{icr}A$,
$\operatorname*{aff}A$, $\operatorname*{cone}A$, and
$\operatorname*{lin}A$ for the closure, the intrinsic core, the
affine hull, the conic hull, and the linear hull of $A$,
respectively. Moreover, by $\overline{\operatorname*{aff}}A,$
$\overline{\operatorname*{cone}}A,$ and
$\overline{\operatorname*{lin}}A$ we denote the closure of the sets
$\operatorname*{aff}A$, $\operatorname*{cone}A$, and
$\operatorname*{lin}A$, respectively. We also use
$\operatorname*{lin}_{0}A$ for the linear space
parallel with the affine hull of $A$, that is, $\operatorname*{lin}%
_{0}A:=\operatorname*{aff}A-a=\operatorname*{lin}(A-a)=\operatorname*{lin}%
(A-A)$ for some (every) $a\in A$, and $\overline{\operatorname*{lin}}%
_{0}A:=\operatorname*{cl}(\operatorname*{lin}_{0}A)$. Clearly $\overline
{\operatorname*{aff}}A=\overline{\operatorname*{aff}}(\operatorname*{cl}A)$
and $\overline{\operatorname*{lin}}_{0}A=\overline{\operatorname*{lin}}%
_{0}(\operatorname*{cl}A)$. For $A,B\subset X$, $a\in X$, $\Gamma
\subset\mathbb{R}$ and $\gamma\in\mathbb{R}$ we set%
\begin{gather*}
A+B:=\{a+b\mid a\in A,\ b\in B\},\quad a+A:=\{a\}+A,\\
\gamma A:=\{\gamma a\mid a\in A\}\text{ for }\gamma\neq0,\quad0A:=\{0\},\quad
\Gamma A:=\cup_{\alpha\in\Gamma}\alpha A,\quad\Gamma a:=\Gamma\{a\}.
\end{gather*}
Therefore, $A+\emptyset=\emptyset+A=\emptyset$ and $0\in\operatorname*{cone}A$
for any set $A\subset X.$

First we recall the notions of quasi interior and quasi relative interior for
convex sets and some properties of these notions.

\medskip
%\subsection{The quasi-interior and quasi relative interior of a convex set}

Let $C\subset X$ be a convex set; the \emph{quasi interior} of $C$ (see
\cite[p.\ 2544]{BorGoe:03}) is the set
\begin{equation*}
\operatorname*{qi}C:=\left\{  x\in C\mid\overline{\operatorname*{cone}%
}(C-x)=X\right\}  ,\label{r-qi}%
\end{equation*}
and the \emph{quasi relative interior} of $C$ (see \cite[Def.\ 2.3]%
{BorLew:92a}) is the set
\[
\operatorname*{qri}C:=\left\{  x\in C\mid\overline{\operatorname*{cone}%
}(C-x)\text{ is a linear space}\right\}  ;
\]
hence $\operatorname*{qi}\emptyset=\operatorname*{qri}\emptyset=\emptyset.$ It
follows (see \cite[Prop.\ 1.2.7 (1.4)]{Zal:02}) that
\[
\operatorname*{qri}C=\left\{  x\in C\mid\overline{\operatorname*{cone}%
}(C-x)=\overline{\operatorname*{lin}}_{0}C\right\}  ;
\]
therefore, because $\overline{\operatorname*{aff}}C=X$ iff $\overline
{\operatorname*{lin}}_{0}C=X,$%
\begin{equation*}
\operatorname*{qi}C=\left\{
\begin{array}
[c]{ll}%
\operatorname*{qri}C & \text{if }\overline{\operatorname*{aff}}C=X,\\
\emptyset & \text{otherwise.}%
\end{array}
\right.  \label{rt9}%
\end{equation*}
Hence (see also \cite[Lem.\ 6]{BotCse:12}),%
\begin{equation}
\operatorname*{qi}C\neq\emptyset\Rightarrow0\in\operatorname*{qi}%
(C-C)\Leftrightarrow\overline{\operatorname*{aff}}C=X\Leftrightarrow
\overline{\operatorname*{lin}}_{0}C=X\Rightarrow\operatorname*{qi}%
C=\operatorname*{qri}C.\label{rt3}%
\end{equation}

Observe (see \cite[Prop.\ 1.2.8 (ii)]{Zal:02}) that for $C\neq\emptyset$ we
have:
\begin{equation}
\lbrack x\in X,\ \text{\ }\overline{\operatorname*{cone}}(C-x)\text{ is a
linear space}]\Rightarrow x\in\operatorname*{cl}C.\label{rt7}%
\end{equation}

Because $\overline{\operatorname*{cone}}(C-x)=\overline{\operatorname*{cone}%
}(\operatorname*{cl}C-x)$ for every $x\in X$, from the definition of
$\operatorname*{qri}C$ and (\ref{rt7}), for $C\neq\emptyset,$ we get%
\begin{equation*}
\operatorname*{qri}(\operatorname*{cl}C)=\left\{  x\in X\mid\overline
{\operatorname*{cone}}(C-x)\text{ is a linear space}\right\}  =\left\{  x\in
X\mid\overline{\operatorname*{cone}}(C-x)=\overline{\operatorname*{lin}}%
_{0}C\right\}  ,\label{rt8}%
\end{equation*}
whence%
\begin{equation}
\operatorname*{qi}C=C\cap\operatorname*{qi}(\operatorname*{cl}C),\quad
\operatorname*{qri}C=C\cap\operatorname*{qri}(\operatorname*{cl}%
C);\label{rt8c}%
\end{equation}
hence, if $A,B\subset X$ are convex sets, then
\begin{equation*}
A\subset B\subset\operatorname*{cl}A\Rightarrow\operatorname*{qri}%
A\subset\operatorname*{qri}B\subset\operatorname*{qri}(\operatorname*{cl}%
A).\label{rt8b}%
\end{equation*}

The facts that $(1-\lambda)C+\lambda\operatorname*{qri}C\subset
\operatorname*{qri}C$ for $\lambda\in(0,1)$ and $\operatorname*{qri}%
(x+C)=x+\operatorname*{qri}C$ for $x\in X$ are well known (see, e.g.,
\cite{BorLew:92a}), and so $\operatorname*{qri}C$ is convex and,
$\operatorname*{cl}(\operatorname*{qri}C)=\operatorname*{cl}C$ if
$\operatorname*{qri}C\neq\emptyset$; therefore, $\overline{\operatorname*{lin}%
}_{0}(\operatorname*{qri}C)=\overline{\operatorname*{lin}}_{0}C$ provided
$\operatorname*{qri}C\neq\emptyset.$ It follows (see \cite[Prop.\ 2.5
(vii)]{BotCseWan:08}) that
\begin{equation}
\operatorname*{qi}(\operatorname*{qi}C)=\operatorname*{qi}C,\quad
\operatorname*{qri}(\operatorname*{qri}C)=\operatorname*{qri}C.\label{rt8d}%
\end{equation}
Indeed, assume that $\operatorname*{qri}C\neq\emptyset.$ By (\ref{rt8c}),
\[
\operatorname*{qri}(\operatorname*{qri}C)=\operatorname*{qri}C\cap
\operatorname*{qri}(\operatorname*{cl}(\operatorname*{qri}%
C))=\operatorname*{qri}C\cap\operatorname*{qri}(\operatorname*{cl}%
C)=\operatorname*{qri}C.
\]

Relation (\ref{rt12a}) in the next result is stated in \cite[Prop.\ 5]%
{GraPop:14} for $D$ a pointed convex cone with $\operatorname*{qi}%
D\neq\emptyset$.

\begin{proposition}
\label{p-gp}Let $C,D\subset X$ be convex sets; then
\begin{align}
C+\operatorname*{qi}D  &  =\operatorname*{qi}(C+\operatorname*{qi}%
D)\subset\operatorname*{qi}(C+D),\label{rt12a}\\
\operatorname*{qri}C+\operatorname*{qri}D  &  =\operatorname*{qri}%
(\operatorname*{qri}C+\operatorname*{qri}D)\subset\operatorname*{qri}(C+D)
\label{rt12b}%
\end{align}

\end{proposition}

Proof. Relation (\ref{rt12a}) is obvious if $\operatorname*{qi}D=\emptyset.$
Assume that $\operatorname*{qi}D\neq\emptyset.$ By (\ref{rt8d}), for $x\in C$
we have that $\operatorname*{qi}(x+\operatorname*{qi}D)=x+\operatorname*{qi}%
(\operatorname*{qi}D)=x+\operatorname*{qi}D\subset\operatorname*{qi}%
(C+\operatorname*{qi}D).$ Hence
\[
\operatorname*{qi}(C+\operatorname*{qi}D)\subset C+\operatorname*{qi}%
D=\cup_{x\in C}(x+\operatorname*{qi}D)\subset\operatorname*{qi}%
(C+\operatorname*{qi}D)\subset\operatorname*{qi}(C+D),
\]
and so (\ref{rt12a}) holds.

Take $x\in\operatorname*{qri}C$ and $y\in\operatorname*{qri}D.$ Since
$\operatorname*{cl}(\operatorname*{qri}C)=\operatorname*{cl}C,$ we have
$\overline{\operatorname*{lin}}_{0}(\operatorname*{qri}C)=\overline
{\operatorname*{lin}}_{0}C$ and $\overline{\operatorname*{cone}}%
(\operatorname*{qri}C-x)=\overline{\operatorname*{lin}}_{0}C;$ similarly for
$D$ and $y.$ Hence%
\begin{align*}
\overline{\operatorname*{cone}}(\operatorname*{qri}C+\operatorname*{qri}%
D)-(x+y)) &  =\operatorname*{cl}[\operatorname*{cone}(\operatorname*{qri}%
C-x)+\operatorname*{cone}(\operatorname*{qri}D-y)]\\
&  =\operatorname*{cl}[\overline{\operatorname*{cone}}(\operatorname*{qri}%
C-x)+\overline{\operatorname*{cone}}(\operatorname*{qri}D-y)]\\
&  =\operatorname*{cl}\left(  \overline{\operatorname*{lin}}_{0}%
C+\overline{\operatorname*{lin}}_{0}D\right)  =\overline{\operatorname*{lin}%
}_{0}(C+D)\\
&  =\overline{\operatorname*{lin}}_{0}(\operatorname*{qri}%
C+\operatorname*{qri}D).
\end{align*}
It follows that $x+y\in\operatorname*{qri}(\operatorname*{qri}%
C+\operatorname*{qri}D).$ Hence $\operatorname*{qri}C+\operatorname*{qri}%
D=\operatorname*{qri}(\operatorname*{qri}C+\operatorname*{qri}D).$ The
inclusion $\operatorname*{qri}C+\operatorname*{qri}D\subset\operatorname*{qri}%
(C+D)$ is well known (see \cite[Lem.\ 3.6 (b)]{BorGoe:03}). Hence
(\ref{rt12b}) holds.  \hfill $\square$

\medskip

Taking into account (\ref{rt3}), from (\ref{rt12a}) one obtains
 \cite[Lem.\ 2.6]{BotCseWan:08} and \cite[Lem.\ 6]{BotCse:12}).

Borwein and Goebel, in \cite[p.\ 2548]{BorGoe:03}, say
\textquotedblleft \textit{Can
$\operatorname*{qri}C+\operatorname*{qri}D$ be a proper subset of
$\operatorname*{qri}(C+D)$? (Almost certainly such sets do
exist.)}\textquotedblright, while Grad and Pop, in \cite[p.\
26]{GraPop:14}, say: \textquotedblleft\textit{we conjecture that in
general when $A,B\subseteq V$ are convex sets with
$\operatorname*{qi}B\neq\emptyset$, it holds
$A+\operatorname*{qi}B=\operatorname*{qi}(A+B)$}\textquotedblright.
The next example answers to both problems mentioned above.

\begin{example}
Take $X:=\ell_{2}:=\big\{(x_{n})_{n\geq1}\subset\mathbb{R}\mid\sum_{n\geq
1}x_{n}^{2}<\infty\big\}$ endowed with its usual norm, $\overline{x}%
:=(n^{-1})_{n\geq1}\in\ell_{2},$ $C:=[0,1]\overline{x}\subset\ell_{2}$ and
$D:=\ell_{1}^{+}:=\big\{(x_{n})_{n\geq1}\subset\mathbb{R}_{+}\mid\sum_{n\geq
1}x_{n}<\infty\big\}\subset\ell_{2}.$ Clearly $C$ and $D$ are convex sets,
$\operatorname*{qri}C=\operatorname*{icr}C=(0,1)\overline{x},$
$\operatorname*{qri}D=\operatorname*{qi}D=\big\{(x_{n})_{n\geq1}\subset
\ell_{1}\mid x_{n}>0\ \forall n\geq1\big\}$ and $\overline{x}\in
\operatorname*{qi}(C+D)=\operatorname*{qri}(C+D),$ but $\overline{x}\notin
C+\operatorname*{qi}D\supset\operatorname*{qri}C+\operatorname*{qri}D.$
\end{example}

Proof. First observe that $\operatorname*{aff}D=\operatorname*{lin}%
_{0}D=D-D=\ell_{1}$ and $\operatorname*{cl}D=\ell_{2}^{+}.$ Therefore,
$\overline{\operatorname*{lin}}_{0}D=\ell_{2},$ and so, using (\ref{rt8c}),%
\begin{equation}
\operatorname*{qri}D=\operatorname*{qi}D=D\cap\operatorname*{qi}\ell_{2}%
^{+}=\left\{  x\in\ell_{1}\mid x_{n}>0~\forall n\geq1\right\}  .\label{rbg2}%
\end{equation}

Clearly, $\ell_{1}^{+}\subset C+D\subset\ell_{2}^{+}+\ell_{2}^{+}=\ell_{2}%
^{+},$ whence $\operatorname*{cl}(C+D)=\ell_{2}^{+}.$ Since $\overline{x}%
\in(C+D)\cap\operatorname*{qi}\ell_{2}^{+},$ we obtain that $\overline{x}%
\in\operatorname*{qi}(C+D)$ using (\ref{rt8c}). Assume that $\overline{x}\in
C+\operatorname*{qi}D.$ Then $\overline{x}\in t\overline{x}+\operatorname*{qi}%
D\subset t\overline{x}+\ell_{1}$ with $t\in\lbrack0,1],$ whence
$(1-t)\overline{x}\in\ell_{1}.$ Because
$\overline{x}\notin\ell_{1},$ we have that $t=1,$ and so
$0\in\operatorname*{qi}D.$ This is a contradiction by (\ref{rbg2}).
The conclusion follows.  \hfill $\square$

\medskip

It is worth observing that for $x_{0}\in C$ we have that%
\begin{equation*}
x_{0}\notin\operatorname*{qi}C\iff\exists x^{\ast}\in X^{\ast}\setminus
\{0\}:\ \inf x^{\ast}(C)=\left\langle x_{0},x^{\ast}\right\rangle ,
\label{rt1b}%
\end{equation*}
that is $x_{0}$ is a support point of $C$, and
\begin{equation}
x_{0}\notin\operatorname*{qri}C\iff\exists x^{\ast}\in X^{\ast}:\ \sup
x^{\ast}(C)>\inf x^{\ast}(C)=\left\langle x_{0},x^{\ast}\right\rangle
\label{rt1}%
\end{equation}
(see also \cite[Lem.\ 2.7]{BorGoe:03}); in particular,
\begin{equation*}
x_{0}\in C\setminus\operatorname*{qri}C\Longrightarrow\exists x^{\ast}\in
X^{\ast}\setminus\{0\}:\ \inf x^{\ast}(C)=\left\langle x_{0},x^{\ast
}\right\rangle . \label{rt2}%
\end{equation*}
Note that in the above implications we do not assume that $\operatorname*{qri}%
C\neq\emptyset.$

\begin{proposition}
\label{p-bg}Let $C,D\subset X$ be nonempty convex sets.

\emph{(i)} If $C\cap\operatorname*{qri}D\neq\emptyset$ then
$\operatorname*{qri}(C\cap D)\subset C\cap\operatorname*{qri}D.$

\emph{(ii)} If $\operatorname*{qri}C\cap\operatorname*{qri}D\neq\emptyset$
then $\operatorname*{qri}(C\cap D)\subset\operatorname*{qri}C\cap
\operatorname*{qri}D.$
\end{proposition}

Proof. (i) Fix $x_{0}\in C\cap\operatorname*{qri}D$ $(\subset C\cap D).$
Consider $x\in\operatorname*{qri}(C\cap D)$ $(\subset C\cap D).$ Suppose that
$x\notin\operatorname*{qri}D.$ By (\ref{rt1}), there exists $x^{\ast}\in
X^{\ast}\setminus\{0\}$ such that
\begin{equation}
\left\langle x,x^{\ast}\right\rangle =\inf x^{\ast}(D)<\sup x^{\ast
}(D).\label{rbg3}%
\end{equation}
Because $x\in C\cap D,$ it follows that $\left\langle x,x^{\ast}\right\rangle
=\inf x^{\ast}(C\cap D).$ Since $x\in\operatorname*{qri}(C\cap D),$ using
again (\ref{rt1}), we have that $\left\langle x,x^{\ast}\right\rangle
=\left\langle y,x^{\ast}\right\rangle $ for every $y\in C\cap D,$ whence
$\left\langle x,x^{\ast}\right\rangle =\left\langle x_{0},x^{\ast
}\right\rangle .$ From (\ref{rbg3}) we get $\left\langle x_{0},x^{\ast
}\right\rangle =\inf x^{\ast}(D)<\sup x^{\ast}(D)$ which implies, by
(\ref{rt1}), that $0\notin\operatorname*{qri}D.$ This contradiction proves
that $\operatorname*{qri}(C\cap D)\subset C\cap\operatorname*{qri}D.$

(ii) Assume that $\operatorname*{qri}C\cap\operatorname*{qri}D\neq\emptyset.$
Then $C\cap\operatorname*{qri}D\neq\emptyset$ and $\operatorname*{qri}C\cap
D\neq\emptyset.$ By (i) we get
\[
\operatorname*{qri}(C\cap D)\subset\left(  C\cap\operatorname*{qri}D\right)
\cap\left(  \operatorname*{qri}C\cap D\right)  =\operatorname*{qri}%
C\cap\operatorname*{qri}D.
\]
The proof is complete.  \hfill $\square$

\medskip

Borwein and Goebel, in \cite[p.\ 2548]{BorGoe:03}, also say
\textquotedblleft\textit{Can $\operatorname*{qri}(C\cap
D)\subset\operatorname*{qri}C\cap
\operatorname*{qri}D$ fail when $\operatorname*{qri}C\cap\operatorname*{qri}%
D\not =\emptyset$?}\textquotedblright\ Proposition \ref{p-bg} (ii)
shows that the answer to this question is negative.


\begin{thebibliography}{1}

\bibitem{BorGoe:03}
Borwein, J., Goebel, R.: Notions of relative interior in {B}anach
spaces.
\newblock J. Math. Sci. (N. Y.) \textbf{115}(4), 2542--2553 (2003)

\bibitem{BorLew:92a}
Borwein, J.M., Lewis, A.S.: Partially finite convex programming.
{I}. {Q}uasi
  relative interiors and duality theory.
\newblock Math. Programming \textbf{57}(1, Ser. B), 15--48 (1992)

\bibitem{BotCse:12}
Bo{\c{t}}, R.I., Csetnek, E.R.: Regularity conditions via
generalized
  interiority notions in convex optimization: new achievements and their
  relation to some classical statements.
\newblock Optimization \textbf{61}(1), 35--65 (2012)

\bibitem{BotCseWan:08}
Bo{\c{t}}, R.I., Csetnek, E.R., Wanka, G.: Regularity conditions via
  quasi-relative interior in convex programming.
\newblock SIAM J. Optim. \textbf{19}(1), 217--233 (2008)

\bibitem{GraPop:14}
Grad, S.M., Pop, E.L.: Vector duality for convex vector optimization
problems
  by means of the quasi-interior of the ordering cone.
\newblock Optimization \textbf{63}(1), 21--37 (2014)

\bibitem{Zal:02}
Z{\u{a}}linescu, C.: Convex analysis in general vector spaces.
\newblock World Scientific Publishing Co. Inc., River Edge, NJ (2002)

\end{thebibliography}
\end{document}